\documentclass[11pt,a4paper]{article}
\usepackage[margin=1in]{geometry}
\usepackage{amsmath,amssymb,amsthm,mathtools}
\usepackage{hyperref}
\usepackage{enumitem}
\usepackage{tikz}
\usepackage{float} 
\usetikzlibrary{arrows.meta,decorations.pathreplacing}
\newtheorem{theorem}{Theorem}[section]
\newtheorem{lemma}[theorem]{Lemma}
\newtheorem{proposition}[theorem]{Proposition}
\newtheorem{remark}[theorem]{Remark}
\newtheorem{example}[theorem]{Example}
\newtheorem{definition}[theorem]{Definition}

\newcommand{\hp}{h_{\mathrm{pol}}}

\newcommand{\NN}{\mathbb N_0}
\def \Dj{\mbox{\raise0.3ex\hbox{-}\kern-0.4em D}}

\title{A flexibility result for polynomial entropy of dendrite maps}
\author{Ma\v{s}a \Dj ori\'c\thanks{Corresponding author: Ma\v{s}a \Dj ori\'c.}\\
Matemati\v{c}ki institut SANU\\
Knez Mihailova 36\\
11000 Beograd\\
Serbia\\
masha@mi.sanu.ac.rs\\ \and
Jelena Kati\'c\\
Matemati\v cki fakultet\\
Studentski trg 16\\
11000 Beograd\\
Serbia\\
jelena.katic@matf.bg.ac.rs
}
\date{}
\begin{document}
\maketitle
\begin{abstract}
For every $\alpha\in(0,1)$ we construct a continuous map on a dendrite whose polynomial entropy is $\alpha$.
\end{abstract}

\medskip


\section{Introduction}
One useful measure of complexity for systems with zero topological entropy is polynomial entropy. The very definitions of topological and polynomial entropy differ in that the former quantifies the overall {\it exponential} growth of orbit complexity, while the latter captures the growth rate on a {\it polynomial} scale. Polynomial entropy serves as a refined tool to distinguish among systems with zero topological entropy. For instance, consider a rotation on the circle versus a homeomorphism with both periodic and wandering points. Although both have topological entropy equal to zero, the first system is evidently simpler than the second. Labrousse proved that polynomial entropy can distinguish between such systems (see Theorem 1 in~\cite{L}).

Topological and polynomial entropy share several properties: both are conjugacy invariants, they depend only on the topology (not on the specific choice of metric), they satisfy the finite union property, and they both obey a product formula. However, there are important differences as well. Certain properties such as the power formula, the $\sigma$-union property, and the variational principle hold for topological entropy but do not necessarily hold for polynomial entropy (see~\cite{M1,M2}). Moreover, topological entropy is equal to the entropy of the system restricted to the non-wandering set, which is a closed, invariant subset. This is not the case for polynomial entropy. Unlike topological entropy, which captures only the behavior on the non-wandering set, polynomial entropy is sensitive to the wandering part of the system.

In \cite{Reg} the authors obtain a rigidity result, as a generalization of a result of Labrousse~\cite{L} for interval and circle homeomorphisms. Namely, we prove that polynomial entropy of a homeomorphism on a dendrite takes only two values - zero and one (see Theorem 13). However, when the dendrite map is only continuous, Roth, Roth and Snoha obtained a flexibility result concerning polynomial entropy. They constructed a dendrite and a surjective continuous map on it, with the given polynomial entropy $\alpha\in (2,3)$ and generalized it for every $\alpha\in[2,+\infty)$(see~\cite{RRS}). That left an open question if there exist a dendrite $X$ and a continuous map $f:X\to X$ such that $\hp(f)\in(0,1)\cup(1,2)$.

We give a (partial) positive answer to this question and for every $\alpha\in(0,1)$ we construct a dendrite $X$ and a continuous map $f$ such that $\hp(f)=\alpha$.

By the Morse–Hedlund theorem, the polynomial entropy of a subshift over a finite alphabet is either zero or at least one. Thus, classical symbolic dynamics cannot provide models with polynomial entropy strictly between zero and one. This obstruction disappears when the alphabet is allowed to be a compact interval. In our construction, we exploit this additional flexibility by considering sequences over the alphabet [0,1] and realizing their dynamics on a dendrite.

\section{Preliminaries}
\subsection{Polynomial entropy}\label{subsec:entropy}

Suppose that $X$ is a compact metric space, and $f:X\rightarrow X$ is continuous. Denote by $d_n^f(x,y)$ the dynamic metric (induced by $f$ and $d$):
$$
d_n^f(x,y)=\max\limits_{0\leqslant k\leqslant n-1}d(f^k(x),f^k(y)).
$$
Fix $Y\subseteq X$. For $\varepsilon>0$, we say that a finite set $E\subset X$ is \textit{$(n,\varepsilon)$-separated} if for every $x,y \in E$, $x\neq y$ it holds $d_n^f(x,y)\geqslant \varepsilon$. Let $\mathrm{sep}(n,\varepsilon;Y)$ denote the maximal cardinality of an $(n,\varepsilon)$-separated set $E$, contained in $Y$.

\begin{definition}\label{def:pol_ent} The \textit{polynomial entropy} of the map $f$ on the set $Y$ is defined by
$$
h_{\mathrm{pol}}(f;Y)=\lim\limits_{\varepsilon \rightarrow 0}\limsup\limits_{n\rightarrow \infty}\frac{\log \mathrm{sep}(n,\varepsilon;Y)}{\log n}.
$$
\end{definition}

We can also define the polynomial entropy as follows. Let $\mathrm{span}(n,\varepsilon;Y)$ denote the minimal number of balls of radius $\varepsilon$ (with respect to $d_n^f$) that cover $Y$ (the centers of the ball are not necessarily elements of $Y$). Denote by $\mathrm{cov}(n,\varepsilon;Y)$ the minimal number of sets $Y_j$ such that the diameters (with respect to $d_n^f$) of $Y_j$ are smaller than $\varepsilon$ and $Y\subseteq\cup_{j=1}^mY_j$.

From the following sequence of inequalities
$$\mathrm{cov}(n,2\varepsilon;Y)\le \mathrm{span}(n,\varepsilon;Y)\le \mathrm{sep}(n,\varepsilon;Y)\le \mathrm{cov}(n,\varepsilon;Y)$$ we conclude that
$$h_{\mathrm{pol}}(f;Y)=\lim\limits_{\varepsilon \rightarrow 0}\limsup\limits_{n\rightarrow \infty}\frac{\log \mathrm{span}(n,\varepsilon;Y)}{\log n}=\lim\limits_{\varepsilon \rightarrow 0}\limsup\limits_{n\rightarrow \infty}\frac{\log \mathrm{cov}(n,\varepsilon;Y)}{\log n}.$$

If $X=Y$ we abbreviate $h_{\mathrm{pol}}(f):=h_{\mathrm{pol}}(f;X)$, $\mathrm{span}(n,\varepsilon):=\mathrm{span}(n,\varepsilon;X)$ etc. We list some of the most relevant properties of the polynomial entropy (for proofs see Propositions $1-4$ in \cite{M2}):
\begin{itemize}
\item[(1)] $h_{\mathrm{pol}}(f^k)=h_{\mathrm{pol}}(f)$, for any $k\geqslant 1$.
\item[(2)] If $Y\subset X$ is a closed, $f$-invariant set, then $h_{\mathrm{pol}}(f;Y)=h_{\mathrm{pol}}(f|_Y)$.
\item[(3)] If $Y=\bigcup_{j=1}^mY_j$ where $Y_j$ are $f$-invariant, then $h_{\mathrm{pol}}(f;Y)=\max\{h_{\mathrm{pol}}(f;Y_j)\mid j=1,\ldots,m\}$\label{finite-union}.\label{(3)}
\item[(4)] If $f:X\to X$, $g:Y\to Y$ and $f\times g:X\times Y\to X\times Y$ is defined as $f\times g (x,y):=(f(x),g(y))$, then $h_{\mathrm{pol}}(f\times g)=h_{\mathrm{pol}}(f)+h_{\mathrm{pol}}(g)$.
\item[(5)] $h_{\mathrm{pol}}(f)$ does not depend on a metric but only on the induced topology.
\item[(6)] $h_{\mathrm{pol}}(\cdot)$ is a \textit{conjugacy invariant} (meaning if $f:X\to X$, $g:X'\to X'$, $\varphi:X\to X'$ is a homeomorphism of compact spaces and $g\circ\varphi=\varphi\circ f$, then $h_{\mathrm{pol}}(f)=h_{\mathrm{pol}}(g)$).
\item[(7)] If $f:X\to X$ and $g:X'\to X'$ are \textit{semi-conjugated}, meaning that $\varphi:X\to X'$ is a continuous surjective map of compact spaces and $g\circ\varphi=\varphi\circ f$, then $h_{\mathrm{pol}}(f)\geqslant h_{\mathrm{pol}}(g)$.
\end{itemize}

\begin{lemma}\label{lem:union}
Let $f:X\to X$ be a continuous map on a compact metric
space, and let $A,B\subseteq X$ be nonempty compact
forward invariant sets. If $B\subseteq\operatorname{Fix}(f)$,
then
\[
h_{\mathrm{pol}}(f|_{A\cup B})
=h_{\mathrm{pol}}(f|_A).
\]
\end{lemma}

\begin{proof}
Since $f|_B=\operatorname{Id}_B$, we have
$h_{\mathrm{pol}}(f|_B)=0$.
By the finite union property of polynomial entropy,
\[
\begin{aligned}
h_{\mathrm{pol}}(f|_{A\cup B})
&=\max\left\{
h_{\mathrm{pol}}(f|_A),
h_{\mathrm{pol}}(f|_B)
\right\}\\
&=h_{\mathrm{pol}}(f|_A).
\end{aligned}
\]
\end{proof}

The following lemma will be useful in our computations.
\begin{lemma}\label{lem:image}
Let $f:X\to X$ be a continuous map on a compact metric space $X$. Then, for every $j\geqslant0$, we have that:
\[
\hp(f)=\hp\bigl(f|_{f^j(X)}\bigr).
\]
\end{lemma}

\begin{proof}
The case $j=0$ is trivial, so assume that $j\geqslant 1$.
Since $f^j(X)\subseteq X$, it follows that $\hp(f|_{f^jX})\leqslant\hp(f)$.

For the reverse inequality, fix $\varepsilon>0$.
Cover $X$ by
\[
C_{j,\varepsilon}
=\operatorname{span}(X,d_j^f,\varepsilon/4)
\]
balls of radius $\varepsilon/4$ in the metric $d_j^f$,
and cover $f^j(X)$ by
\[
\operatorname{span}(f^j(X),d_n^f,\varepsilon/4)
\]
balls of radius $\varepsilon/4$ in the metric $d_n^f$.

For each nonempty intersection of an element of the
first cover with the preimage under $f^j$ of an
element of the second cover, choose a representative
point of $X$.

Any two points $x,y$ belonging to such an intersection
satisfy
\[
d_j^f(x,y)<\frac{\varepsilon}{2}
\quad\text{and}\quad
d_n^f(f^j(x),f^j(y))<\frac{\varepsilon}{2}.
\]
Since
\[
d_{n+j}^f(x,y)
=\max\{d_j^f(x,y),d_n^f(f^j(x),f^j(y))\},
\]
we obtain
\[
d_{n+j}^f(x,y)<\varepsilon.
\]

Therefore, the Bowen balls of radius $\varepsilon$
centered at the chosen representative points cover
$X$, and hence
\begin{equation}\label{eq:iterated-image-cover}
\operatorname{span}(X,d_{n+j}^f,\varepsilon)
\leqslant C_{j,\varepsilon}
\operatorname{span}(f^j(X),d_n^f,\varepsilon/4).
\end{equation}

Since $j$ and $C_{j,\varepsilon}$ are independent
of $n$, taking logarithms, dividing by $\log n$,
and passing to the limit superior as $n\to\infty$,
we obtain
\[
\begin{aligned}
&\limsup_{n\to\infty}
\frac{\log\operatorname{span}(X,d_n^f,\varepsilon)}
{\log n}\\
&\qquad\leqslant
\limsup_{n\to\infty}
\frac{\log\operatorname{span}
(f^j(X),d_n^f,\varepsilon/4)}
{\log n}.
\end{aligned}
\]
Finally, letting $\varepsilon\to 0$, we conclude that
\[
h_{\mathrm{pol}}(f)
\leqslant h_{\mathrm{pol}}(f|_{f^j(X)}),
\]
which completes the proof.
\end{proof}

\subsection{Shift spaces and polynomial entropy}

Unlike in classical symbolic dynamics, where the alphabet is typically finite and discrete, we allow the alphabet to be an arbitrary compact metric space. Let $(A,d_A)$ be a compact metric space.
Without loss of generality, we may assume that
$\operatorname{diam}(A)\leqslant 1$.
Consider the space
\[
\Sigma_A=A^{\mathbb N_0}
\]
endowed with the metric
\[
d(x,y)=\sup_{r\geqslant 0}2^{-r}d_A(x_r,y_r),
\]
which induces the product topology.

The shift map $\sigma\colon\Sigma_A\to\Sigma_A$
is defined by
\[
\sigma(x)_r=x_{r+1},\qquad r\geqslant 0.
\]

For each $n\geqslant 1$, define the pseudometric
$\rho_n$ on $\Sigma_A$ by
\[
\rho_n(x,y)=
\max_{0\leqslant r<n}d_A(x_r,y_r).
\]

Thus, $\rho_n$ measures the maximal distance
between the corresponding coordinates of two
sequences during the first $n$ steps.

For a subset $Z\subseteq\Sigma_A$, we denote by
\[
\operatorname{span}(Z,\rho_n,\varepsilon)
\]
the minimal number of $\rho_n$-balls of radius
$\varepsilon$, with centers in $\Sigma_A$,
needed to cover $Z$.

The following lemma shows that polynomial entropy
can be computed using the pseudometrics $\rho_n$
instead of the Bowen metrics associated with $d$.

\begin{lemma}\label{lem:productmetric}
Let $Z\subseteq\Sigma_A$ be a nonempty compact
$\sigma$-invariant set. Then
\[
\hp(\sigma|_Z)
=
\lim_{\varepsilon\to 0}
\limsup_{n\to\infty}
\frac{
\log\operatorname{span}(Z,\rho_n,\varepsilon)
}{
\log n
}.
\]
\end{lemma}

\begin{proof}
Let
\[
d_n^\sigma(x,y)
=
\max_{0\leqslant i<n}d(\sigma^i(x),\sigma^i(y))
\]
be the Bowen metric associated with $d$.

By the definition of $d$, we have
\[
\rho_n(x,y)\leqslant d_n^\sigma(x,y).
\]
Consequently,
\[
\operatorname{span}(Z,\rho_n,\varepsilon)
\leqslant
\operatorname{span}(Z,d_n^\sigma,\varepsilon).
\]

Conversely, fix $\varepsilon>0$ and choose
$m\geqslant 1$ such that $2^{-m}<\varepsilon$.
For every $0\leqslant i<n$, we have
\[
\begin{aligned}
d(\sigma^i(x),\sigma^i(y))
&=
\sup_{r\geqslant 0}2^{-r}
d_A(x_{i+r},y_{i+r})\\
&\leqslant
\max\{\rho_{n+m}(x,y),2^{-m}\}.
\end{aligned}
\]
Hence
\[
d_n^\sigma(x,y)
\leqslant
\max\{\rho_{n+m}(x,y),2^{-m}\}.
\]

It follows that
\[
\operatorname{span}(Z,d_n^\sigma,\varepsilon)
\leqslant
\operatorname{span}(Z,\rho_{n+m},\varepsilon).
\]

Combining the two inequalities, we obtain
\[
\begin{aligned}
\operatorname{span}(Z,\rho_n,\varepsilon)
&\leqslant
\operatorname{span}(Z,d_n^\sigma,\varepsilon)\\
&\leqslant
\operatorname{span}(Z,\rho_{n+m},\varepsilon).
\end{aligned}
\]

Since $m$ depends only on $\varepsilon$ and not
on $n$, and
\[
\lim_{n\to\infty}
\frac{\log(n+m)}{\log n}=1,
\]
the three quantities have the same polynomial
growth rate.

Finally, allowing the centers of spanning balls
to belong to $\Sigma_A$ rather than to $Z$
does not change the polynomial entropy, since
an ambient spanning ball intersecting $Z$ can
be replaced by a ball centered in $Z$ with
twice the radius.

Letting $\varepsilon\to 0$ completes the proof.
\end{proof}

In particular, if $A$ is a finite alphabet endowed
with the discrete metric, then
\[
\rho_n(x,y)=
\max_{0\leqslant r<n}d_{\mathrm{disc}}(x_r,y_r).
\]
If $A=[0,1]$ is endowed with the usual metric, then
\[
\rho_n(x,y)=
\max_{0\leqslant r<n}|x_r-y_r|.
\]
We will use the latter pseudometrics when studying
horizontal histories associated with the dendrite
construction.

\begin{remark}\label{rem:closure}
Let $G$ be a dense subset of $Z\subseteq\Sigma_A$.
For every $n\geqslant 1$ and every
$0<\varepsilon<\varepsilon'$, we have
\[
\operatorname{span}(Z,\rho_n,\varepsilon')
\leqslant
\operatorname{span}(G,\rho_n,\varepsilon).
\]

Indeed, let $m=\operatorname{span}(G,\rho_n,\varepsilon)$ and choose points $w_1,\ldots,w_m\in\Sigma_A$
such that
\[
G\subseteq
\bigcup_{i=1}^{m}B_{\rho_n}(w_i,\varepsilon).
\]
Since $\rho_n$ is continuous, the sets $C_i=\{z\in\Sigma_A:\rho_n(z,w_i)\leq\varepsilon\}$ are closed. Since their union is finite and
contains $G$, it also contains $\overline{G}$. Therefore,
\[
Z\subseteq\overline{G}
\subseteq\bigcup_{i=1}^{m}C_i
\subseteq
\bigcup_{i=1}^{m}B_{\rho_n}(w_i,\varepsilon'),
\]
where the last inclusion follows from
$\varepsilon<\varepsilon'$. Thus, the same centers $w_1,\ldots,w_m$,
with the radius increased from $\varepsilon$
to $\varepsilon'$, provide a cover of $Z$. This proves the desired inequality.

We will use this observation several times
throughout the paper. Namely, this elementary observation justifies the radius enlargements in Lemmas \ref{lem:upper}, \ref{lem:Qcover}, and \ref{lem:perturb}. In particular, no factor depending on $n$ is introduced.
\end{remark}

\section{Main results}

\subsection{Construction of dendrite and dendrite map}
Let $(a_j)_{j\ge0}$ be pairwise distinct numbers in $(0,1)$ satisfying
\begin{equation}\label{eq:increments}
 |a_{j+1}-a_j|\longrightarrow0,
\end{equation}
and let $h_j>0$ with $h_j\to0$. Put
\[
K=[0,1]\times\{0\},\qquad q_j=(a_j,0),\qquad p_j=(a_j,h_j),\qquad
I_j=\{a_j\}\times[0,h_j],
\]
and
\begin{equation}\label{eq:dendrite}
D=K\cup\bigcup_{j\ge0}I_j,\qquad
D_j=K\cup\bigcup_{i\ge j}I_i\quad(j\ge0).
\end{equation}
We use the Euclidean metric restricted from $\mathbb R^2$. The notation $[u,v]_{D}$ stands for the unique arc joining $u$ to $v$ in $D$. We call $I_j$ a {\it hair}.

\begin{proposition}\label{prop:dendrite}
$D$ is a dendrite and each $D_j$ is a subdendrite.
\end{proposition}
\begin{proof}

$D$ is bounded. Let us prove that it is also closed. Let $x_m\in D$ be a sequence of points in $D$. If infinitely many $x_m$ belong to $K$ or to a fixed finite union of hairs (which are both compact), there exists a convergent subsequence with a limit in D. Otherwise, take a subsequence, $x_m\in I_{j_m}$ with $j_m\to\infty$, and its convergent subsequence, denote it again by $x_m$. Its vertical coordinate is between $0$ and $h_{j_m}$ and hence tends to zero. Its horizontal coordinate belongs to $[0,1]$, so the limit belongs to $K$. Again, we conclude there is a convergent subsequence with a limit in $D$. $D$ is bounded and closed, and therefore compact. Also, $D$ is the union of connected hairs, each intersecting the connected set $K$, so $D$ is also connected.

Any two points are joined by a unique arc: it consists of the appropriate vertical segments and the appropriate subinterval of $K$. Consequently, $D$ contains no simple closed curve.

Let $(a_j,t)\in D\backslash K$, $t>0$. Since $h_j\to 0$, only finitely many hairs exceed the height $\frac{t}{2}$. Since $a_i\neq a_j$ for $i\neq j$, there exists $\delta>0$ such that for all $0<r<\delta$ we have that $B(z,r)\cap D=B(z,r)\cap I_j$, which is a connected set. Now, let $(a,0)\in K$. For any $\varepsilon>0$, the set $U_{\varepsilon}=D\cap ((a-\varepsilon,a+\varepsilon)\times(-\varepsilon,\varepsilon))$ is connected. We can now conclude that $D$ is locally connected. Since $D$ is a locally connected continuum which contains no simple closed curves, $D$ is indeed a dendrite.

In the same way, we prove that $D_j$ is a subdendrite.
\end{proof}

For each $j$, choose a homeomorphism
\begin{equation}\label{eq:hairmap}
\varphi_j:I_j\rightarrow [q_j,p_{j+1}]_{D}
=[q_j,q_{j+1}]_K\cup I_{j+1},
\end{equation}
with $\varphi_j(q_j)=q_j$ and $\varphi_j(p_j)=p_{j+1}$. Define
\begin{equation}\label{eq:F}
F|_K=\mathrm{Id}_K,\qquad F|_{I_j}=\varphi_j.
\end{equation}
The definitions agree at $q_j$. 

\begin{figure}[H]
    \centering 
    \includegraphics[width=0.9\textwidth]{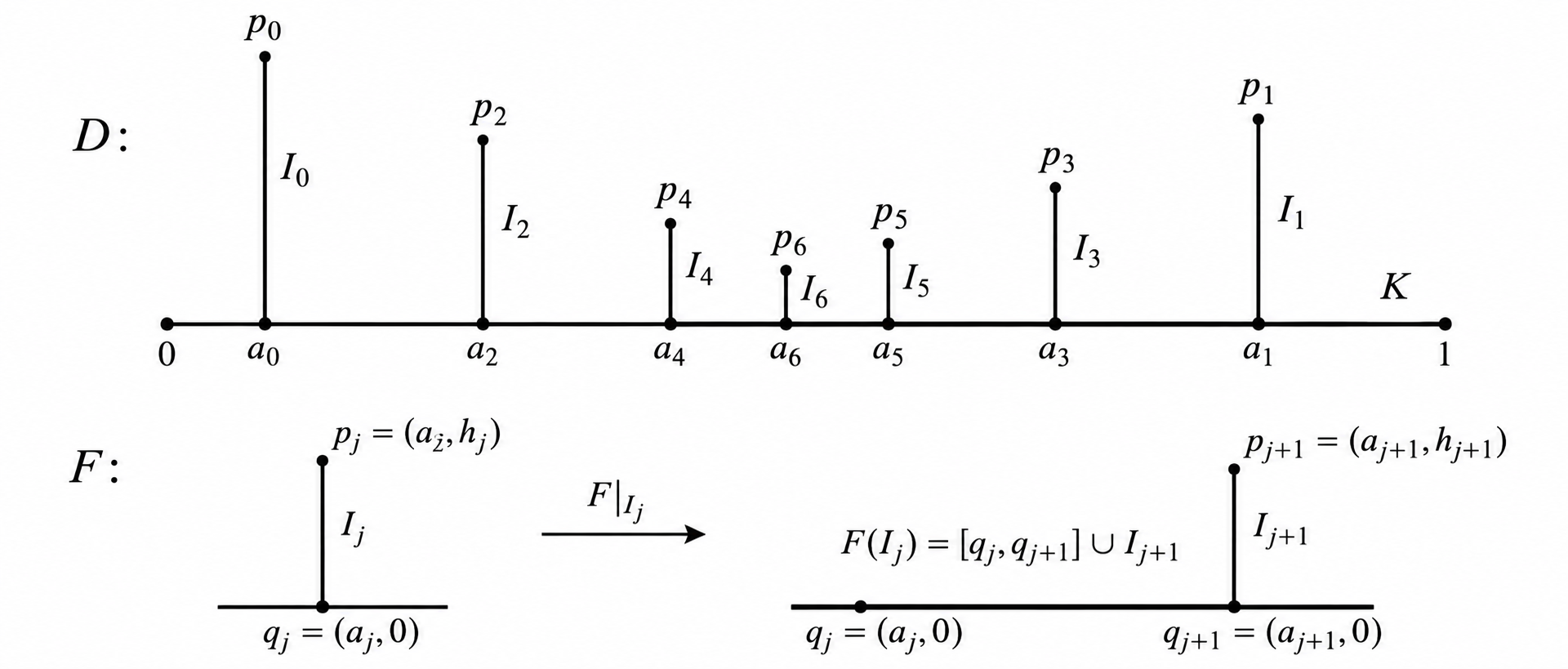}
\end{figure}

\begin{lemma}\label{lem:continuity}
$F:D\to D$ is a continuous map. $F(D_j)=D_{j+1}$, for all $j\geqslant0$ and hence $F^j(D)=D_j$.

\end{lemma}
\begin{proof}

For every $j\geqslant 0$ and every $x\in I_j$, we have
\[
d(x,q_j)\leq h_j
\]
and, since $F(I_j)=[q_j,q_{j+1}]\cup I_{j+1}$,
\[
d(F(x),q_j)\leq |a_{j+1}-a_j|+h_{j+1}.
\]
Hence, by the triangle inequality,
\[
\sup_{x\in I_j}d(F(x),x)
\leq h_j+|a_{j+1}-a_j|+h_{j+1}.
\]
Since $h_j\to 0$ and $|a_{j+1}-a_j|\to 0$, it follows that
\[
\lim_{j\to\infty}\sup_{x\in I_j}d(F(x),x)=0.
\]

For each $N\in\mathbb{N}$, define $F_N:D\to D$ by
\[
F_N(x)=
\begin{cases}
F(x), & x\in K\cup\displaystyle\bigcup_{j=0}^{N}I_j,\\
x, & \text{otherwise}.
\end{cases}
\]
The sets
\[
A_N=K\cup\bigcup_{j=0}^{N}I_j
\quad\text{and}\quad
B_N=K\cup\bigcup_{j>N}I_j
\]
are compact. Indeed, $A_N$ is a finite union of closed subsets of $D$, and the proof that $B_N$ is compact is the same as for $D$. It holds $A_N\cup B_N=D$ and 
$A_N\cap B_N=K$.

The restriction $F_N|_{A_N}$ is continuous by the 
pasting lemma, since it is obtained by gluing finitely 
many continuous maps that agree on their intersections.
Moreover, $F_N|_{B_N}=\operatorname{Id}_{B_N}$.
Since these two restrictions agree on $K$, another 
application of the pasting lemma yields that $F_N$ 
is continuous on $D$.

Furthermore,
\[
\begin{aligned}
\sup_{x\in D}d(F_N(x),F(x))
&=\sup_{j>N}\sup_{x\in I_j}d(x,F(x))\\
&\leq\sup_{j>N}
\left(h_j+|a_{j+1}-a_j|+h_{j+1}\right)\\
&\longrightarrow 0
\end{aligned}
\]
as $N\to\infty$. Therefore, $F$ is the uniform limit of the continuous 
maps $F_N$, and hence it is continuous.

Since $F(K)=K$ and $F(I_j)=[q_j,q_{j+1}]_K\cup I_{j+1}$,
\[
F(D_j)=K\cup\bigcup_{i\geqslant j}I_{i+1}=D_{j+1}.
\]
\end{proof}

\subsection{History space}
Let $\Sigma=[0,1]^{\NN}$, equipped with its product topology, and let $\sigma:\Sigma\to\Sigma$ be the left shift. Write $c^\infty=(c,c,c,\ldots)$. Define the horizontal projection $\pi:D\to[0,1]$ by $\pi(x,y)=x$, and the history map
\begin{equation}\label{eq:Psi}
\Psi:D\to\Sigma,\qquad
\Psi(x)=\bigl(\pi(x),\pi(F x),\pi(F^2x),\ldots\bigr).
\end{equation}
Each coordinate of $\Psi$ is continuous, so $\Psi$ is continuous in the product topology. Moreover,
\begin{equation}\label{eq:factor}
\Psi\circ F=\sigma\circ\Psi.
\end{equation}
Define the \emph{full horizontal-history space}
\begin{equation}\label{eq:X=psi(D)}
    X=\Psi(D).
\end{equation}
It is compact and shift invariant in the forward sense, $\sigma(X)\subseteq X$.

For $0\le j\le s<\infty$, define the \emph{stopped history}
$$u^{j,s}=(a_j,a_{j+1},\ldots,a_s,a_s,a_s,\ldots),$$
and the \emph{unstopped history}
$$u^{j,\infty}=(a_j,a_{j+1},a_{j+2},\ldots).$$
The ordinary-history space is defined by
\[Z=\overline{\{u^{j,\infty}:j\geqslant0\}}\]
The \emph{stopped-history space} is
\begin{equation}\label{eq:Y}
Y=\overline{\{u^{j,s}:0\le j\le s\le\infty\}}\ \subseteq\Sigma,
\end{equation}
where the closure is in the product topology. It is compact and forward shift invariant. Clearly, $Z\subseteq Y$. 

Since $|a_{j+1}-a_j|\to 0$, the space $Y$ consists of the following three types of sequences:
\begin{enumerate}
    \item[(i)] stopped histories $u^{j,s}$, where 
    $0\leqslant j< s<\infty$;
    
    \item[(ii)] unstopped histories $u^{j,\infty}$, 
    where $j\geqslant 0$;
    
    \item[(iii)] constant sequences $(c,c,\ldots)$, 
    where $c\in\overline{\{a_j:j\geqslant 0\}}$.
\end{enumerate}

Indeed, let $(u^{j_n,s_n})_{n\geqslant 0}$ be a convergent 
sequence of generators. Passing to a subsequence, 
we may assume that either $j_n\to\infty$ or $j_n=j$ 
for some fixed $j\geqslant 0$.

In the former case, for every $r\geqslant 0$, we have
$$\left|u^{j_n,s_n}_{r+1}-u^{j_n,s_n}_r\right|
\leqslant \sup_{i\geqslant j_n}|a_{i+1}-a_i|.$$
Since $j_n\to\infty$ and $|a_{i+1}-a_i|\to0$ the right-hand side tends to zero. As convergence in the product topology implies coordinatewise convergence, any two consecutive coordinates of the limit sequence coincide. Hence the limit is a constant sequence $c^\infty$, where $c\in\overline{\{a_j:j\geqslant0\}}$.

In the latter case, passing to a further subsequence, 
we may assume that either $s_n=s<\infty$ for every $n$ 
or $s_n\to\infty$. The corresponding limit is 
$u^{j,s}$ or $u^{j,\infty}$, respectively.

Conversely, every constant sequence $(c,c,\ldots)$ 
with $c\in\overline{\{a_j:j\geqslant 0\}}$ is a limit 
of generators of the form $u^{j,j}$.

\begin{proposition}\label{prop:realization}
Every generator $u^{j,s}$, $j\le s\le\infty$, belongs to $X$. Consequently $Y\subseteq X$.
\end{proposition}
\begin{proof}
For $s=j$, take $x=q_j$, which is fixed. For finite $s>j$, put $x_s=q_s$ and recursively define
\[
x_r=\bigl(F|_{I_r}\bigr)^{-1}(x_{r+1}),\qquad r=s-1,\ldots,j.
\]
The inverse is well defined because $x_{r+1}\in I_{r+1}\subseteq F(I_r)$. The orbit of $x_j$ follows $x_j,x_{j+1},\ldots,x_s=q_s$ and then remains at $q_s$, hence $\Psi(x_j)=u^{j,s}$. For $s=\infty$, use $x=p_j$: by construction $F^r(p_j)=p_{j+r}$ for every $r$, giving $\Psi(p_j)=u^{j,\infty}$. Since $X$ is closed, it contains the closure of these histories.
\end{proof}

\begin{remark}
The inclusion $Y\subseteq X$ is always strict.
Indeed, since $a_0\neq a_1$, choose $c$ strictly
between $a_0$ and $a_1$. By construction,
\[
(c,0)\in[q_0,q_1]\subseteq F(I_0).
\]
Hence there exists $x\in I_0$ such that
$F(x)=(c,0)$. Since $(c,0)$ is a fixed point,
we obtain
\[
\Psi(x)=(a_0,c,c,\ldots)\in X.
\]
On the other hand, using the description of $Y$ above
and the fact that $a_i\neq a_j$, for $i\neq j$,
every element of $Y$ whose first coordinate is $a_0$
has its second coordinate equal to either $a_0$ or $a_1$.
Since $c\notin\{a_0,a_1\}$, we conclude that
\[
(a_0,c,c,\ldots)\notin Y.
\]
Therefore, $Y\subsetneq X$.
\end{remark}

\subsection{Main theorem}

We are now ready to prove the main result of this section. Our aim is to show that the polynomial entropy of $F$ can be computed entirely in terms of stopped histories. We first show that passing from the dendrite $D$ to the full horizontal-history space $X$ does not change the polynomial entropy. We then prove that, although $Y$ is a proper subset of $X$, the two spaces have the same polynomial entropy under the shift map.

\begin{lemma}\label{lem:factor-equality}
$\hp(F)=\hp(\sigma|_{X})$.
\end{lemma}

\begin{proof}
Since $\Psi:D\to X=\Psi(D)$ is a continuous factor map, we have
\[
h_{\mathrm{pol}}(\sigma|_X)
\leq h_{\mathrm{pol}}(F).
\]
Let us prove the reverse inequality. Put $H_j=\sup_{i\geqslant j}h_i\to0$. Let $x,y\in D_j$. Since $D_j$ is forward invariant,
the vertical coordinates of $F^k(x)$ and $F^k(y)$
belong to $[0,H_j]$ for every $k\geqslant 0$.
Using the Euclidean metric inherited from $\mathbb{R}^2$,
we obtain
\[
\|F^k(x)-F^k(y)\|
\leq
|\pi(F^k(x))-\pi(F^k(y))|+H_j.
\]
Consequently,
\begin{equation}\label{eq:bowen-comparison}
d_n^F(x,y)
\leq \rho_n(\Psi(x),\Psi(y))+H_j.
\end{equation}

Fix $\varepsilon>0$ and choose $j$ sufficiently large
so that $H_j<\frac{\varepsilon}{2}$. Let
\[
m=\operatorname{span}
\left(\Psi(D_j),\rho_n,\frac{\varepsilon}{4}\right).
\]
By choosing a point of $\Psi(D_j)$ in each nonempty covering ball and doubling the radius, we obtain
points $u_1,\ldots,u_m\in\Psi(D_j)$ such that
\[
\Psi(D_j)\subseteq
\bigcup_{i=1}^{m}
B_{\rho_n}\left(u_i,\frac{\varepsilon}{2}\right).
\]

For each $i\in\{1,\ldots,m\}$, choose a point
$x_i\in D_j$ satisfying
\[
\Psi(x_i)=u_i.
\]

We claim that the Bowen balls
$B_{d_n^F}(x_i,\varepsilon)$, $i=1,\ldots,m$,
cover $D_j$.
Indeed, for any $x\in D_j$, there exists an index
$i$ such that
\[
\rho_n(\Psi(x),\Psi(x_i))
=\rho_n(\Psi(x),u_i)
<\frac{\varepsilon}{2}.
\]
By \eqref{eq:bowen-comparison}, we obtain
\[
\begin{aligned}
d_n^F(x,x_i)
&\leq \rho_n(\Psi(x),\Psi(x_i))+H_j\\
&<\frac{\varepsilon}{2}
+\frac{\varepsilon}{2}\\
&<\varepsilon.
\end{aligned}
\]
Therefore,
\begin{equation}\label{eq:tail-cover}
\operatorname{span}(D_j,d_n^F,\varepsilon)
\leq
\operatorname{span}
\left(\Psi(D_j),\rho_n,\frac{\varepsilon}{4}\right).
\end{equation}

Since $D_j=F^j(D)$, using the estimate~(\ref{eq:iterated-image-cover}) obtained
in the proof of Lemma \ref{lem:image} and inequality \eqref{eq:tail-cover}, we get
\[
\begin{aligned}
\operatorname{span}(D,d_{n+j}^F,4\varepsilon)
&\leq C_{j,4\varepsilon}
\operatorname{span}(D_j,d_n^F,\varepsilon)\\
&\leq C_{j,4\varepsilon}
\operatorname{span}
\left(\Psi(D_j),\rho_n,\frac{\varepsilon}{4}\right)\\
&\leq C_{j,4\varepsilon}
\operatorname{span}
\left(X,\rho_n,\frac{\varepsilon}{4}\right).
\end{aligned}
\]
Here, the last inequality follows from $\Psi(D_j)\subseteq X$ and the fact that ambient covering balls are allowed.
Since $j$ and $C_{j,4\varepsilon}$ are independent
of $n$, taking logarithms, dividing by $\log n$,
passing to the limit superior as $n\to\infty$,
and then letting $\varepsilon\to 0$, we obtain
\[
h_{\mathrm{pol}}(F)
\leq h_{\mathrm{pol}}(\sigma|_X).
\]
\end{proof}

Let us define the following:
$$
\Delta_j=\sup_{i\geqslant j}|a_{i+1}-a_i|\rightarrow0.
$$
\begin{lemma}\label{lem:approx}
For every $j\geqslant0$ and every $x\in D_j\backslash K$, there exists $v\in Y$ such that
\begin{equation}\label{eq:approx}
\sup_{r\geqslant0}|\Psi(x)_r-v_r|\leqslant\Delta_j.
\end{equation}
\end{lemma}
\begin{proof}
Suppose $x\in I_i\setminus K$ with $i\geqslant j$. At each step, a point off the base in $I_i$ either maps into $I_{i+1}\setminus K$ or enters $K$; once in $K$, its orbit is constant.

If the orbit never enters $K$, then $\Psi(x)=u^{i,\infty}\in Y$, so we can take $v=u^{i,\infty}$. Otherwise, let $m\geqslant1$ be the first entrance time. The horizontal history has the form
\[
\Psi(x)=(a_i,a_{i+1},\ldots,a_{i+m-1},c,c,\ldots),
\]
where $c$ belongs to the closed interval with endpoints $a_{i+m-1}$ and $a_{i+m}$, because the last image lands on the base segment of $F(I_{i+m-1})$. Take
\[
v=u^{i,i+m}=(a_i,\ldots,a_{i+m-1},a_{i+m},a_{i+m},\ldots)\in Y.
\]
The histories agree before time $m$; at every later time their discrepancy is
$|c-a_{i+m}|\leqslant|a_{i+m}-a_{i+m-1}|\leqslant\Delta_j$. This proves \eqref{eq:approx}.
\end{proof}

\begin{lemma}\label{lem:XY}
$\hp(\sigma|_{X})=\hp(\sigma|_{Y})$.
\end{lemma}

\begin{proof}
Since $Y\subseteq X$, by monotonicity we have
\[
\hp(\sigma|_{Y})
\leqslant \hp(\sigma|_{X}).
\]
For the reverse inequality, fix $\varepsilon>0$
and choose $j$ sufficiently large so that
\[
\Delta_j<\frac{\varepsilon}{4}.
\]
By Lemma \ref{lem:approx}, for every $u\in\Psi(D_j\backslash K)$ there
exists $v\in Y$ such that for every $n\geqslant 1$:
\[
\rho_n(u,v)\leqslant\Delta_j
<\frac{\varepsilon}{4}.
\]
Let $m=\operatorname{span}
\left(Y,\rho_n,\frac{\varepsilon}{4}\right).$ There exist ambient centers $v_1,\ldots,v_m\in\Sigma$
such that
\[
Y \subseteq
\bigcup_{i=1}^{m}
B_{\rho_n}\left(v_i,\frac{\varepsilon}{4}\right).
\]
For any $u\in\Psi(D_j\backslash K)$, choose $v\in Y$
as above and an index $i$ such that $\rho_n(v,v_i)<\frac{\varepsilon}{4}.$ By the triangle inequality, we obtain
\[
\rho_n(u,v_i)
\leqslant \rho_n(u,v)+\rho_n(v,v_i)
<\frac{\varepsilon}{2}.
\]
Consequently,
\[
\Psi(D_j\backslash K)\subseteq
\bigcup_{i=1}^{m}
B_{\rho_n}\left(v_i,\frac{\varepsilon}{2}\right).
\]
and furthermore
\begin{equation}\label{eq:stopped-history-cover}
\operatorname{span}(\Psi(D_j\backslash K),\rho_n,\varepsilon)
\leqslant
\operatorname{span}
\left(Y,\rho_n,\frac{\varepsilon}{4}\right).
\end{equation}

Since $\Psi(K)$ consists of constant sequences
and
\[
\rho_n(c^\infty,d^\infty)=|c-d|,
\]
there exists a constant $M_\varepsilon$,
independent of $n$, such that
\[
\operatorname{span}(\Psi(K),\rho_n,\varepsilon)
\leqslant M_\varepsilon.
\]

Therefore,
\[
\operatorname{span}(\Psi(D_j),\rho_n,\varepsilon)
\leqslant
\operatorname{span}(Y,\rho_n,\varepsilon/4)
+M_\varepsilon.
\]

Since $\Psi(D_j)=\sigma^j(X)$, the same
covering argument as in the proof of Lemma \ref{lem:image} gives
\[
\begin{aligned}
\operatorname{span}(X,\rho_{n+j},4\varepsilon)
&\leqslant C_{j,4\varepsilon}
\operatorname{span}(\Psi(D_j),\rho_n,\varepsilon)\\
&\leqslant C_{j,4\varepsilon}(
\operatorname{span}
\left(Y,\rho_n,\frac{\varepsilon}{4}\right)+M_{\varepsilon}),
\end{aligned}
\]
where $C_{j,4\varepsilon}$ and $M_{\varepsilon}$ are independent of $n$.

Taking logarithms, dividing by $\log n$,
passing to the limit superior as $n\to\infty$,
and then letting $\varepsilon\to 0$, we obtain
\[
h_{\mathrm{pol}}(\sigma|_{X})
\leqslant h_{\mathrm{pol}}(\sigma|_{Y}).
\]
\end{proof}

\begin{theorem}\label{thm:realization}
Let $(a_j)_{j\geqslant0}$ be pairwise distinct in $(0,1)$ with $|a_{j+1}-a_j|\to0$, and let $h_j>0$ satisfy $h_j\to0$. Let $D$ be the dendrite defined in \eqref{eq:dendrite}, $F:D\to D$ the map defined in \eqref{eq:F}, and $X$ and $Y$ the associated horizontal-history spaces as in \eqref{eq:X=psi(D)} and \eqref{eq:Y}. Then
\begin{equation}\label{eq:main}
{\hp(F)=\hp(\sigma|_{X})=\hp(\sigma|_{Y}).}
\end{equation}
In particular, the value is independent of the positive hair heights $h_j\to0$ and of the particular homeomorphisms $\varphi_j$ satisfying the specified endpoint conditions.
\end{theorem}
\begin{proof}
The first equality is Lemma~\ref{lem:factor-equality} and the second equality is Lemma~\ref{lem:XY}.
\end{proof}

Let us now show a simple application of our previous construction, namely the one where $(a_j)_{j\geqslant0}$ is a convergent sequence.

\begin{example}
Suppose $a_j\to c$. We will show that $\hp(\sigma|_{Y})=0$ and hence $\hp(F)=0$.  Fix $\varepsilon>0$ and choose $j$ so that $|a_k-c|<\varepsilon/4$ for all $k\geqslant j$. For every generator $u^{i,s}$ with $i\geqslant j$, we have
$$\rho_n(u^{i,s},c^\infty)
\leqslant\sup_{k\geqslant j}|a_k-c|
\leqslant\frac{\varepsilon}{4}$$
for every $n\geqslant1$.

For each of the finitely many $i<j$, consider the generators $u^{i,s}$ with $j\leqslant s<\infty$. Since the first $s-i+1$ coordinates of $u^{i,s}$ and $u^{i,\infty}$ agree and all subsequent coordinates belong to the $\varepsilon$-neighbourhood of $c$, we have
$$\rho_n(u^{i,s},u^{i,\infty})
\leq\frac{\varepsilon}{2}$$
for every $n\geqslant1$.

The remaining generators $u^{i,s}$, with $i\leqslant s<j$, form a finite set. Consequently, the entire generating family can be covered by \(\rho_n\)-balls of radius $\varepsilon/2$ (taking closed balls), centred at the points
$$c^\infty,\quad
u^{i,\infty}\ (0\leqslant i<j),\quad
u^{i,s}\ (0\leqslant i\leq s<j).$$

Since $Y$ is the closure of the generating family and $\rho_n$ is continuous, we obtain
$$\operatorname{span}(Y,\rho_n,\varepsilon)
\leqslant 1+j+\frac{j(j+1)}{2}.$$
Since $j$ depends only on $\varepsilon$ and not on $n$, it follows that $h_{\mathrm{pol}}(\sigma|_Y)=0$.
\end{example}

\subsection{Construction of auxiliary sequence}

We first construct an auxiliary sequence $b=(b_j)_{j\geqslant 0}$
consisting of long constant blocks, alternating between
$0$ and $1$, with linear transitions between them.
The main point of the construction is the choice of
the lengths of these blocks and transitions.
If $N_k$ is the length of the $k$th constant block
and $L_k$ is the length of the following transition,
we require that
\[
    \frac{N_k}{L_k}\asymp N_k^\alpha.
\]
This will allow us to obtain the desired polynomial
entropy.  This deliberate simplification allows exact plateaus, which make the estimates considerably easier.
We will then slightly perturb the sequence
$b$ to obtain a sequence $a$ with pairwise distinct
entries, while preserving the polynomial entropy.

Fix $0<\alpha<1$ and set
\[N_k=2^{2^k},\quad L_k=\lceil N_k^{1-\alpha}\rceil,\quad e_0=0,\quad t_k=e_{k-1}+N_k,\quad e_k=t_k+L_k.\]
Let $c_k=0$ for odd $k$ and $c_k=1$ for even $k$. Define
\[b_j=\begin{cases}c_k,&e_{k-1}\le j\le t_k,\\c_k+(c_{k+1}-c_k)(j-t_k)/L_k,&t_k\le j\le e_k.\end{cases}\]
We call the integer interval $P_k:=[e_{k-1},t_k]\cap\mathbb N_0$
the $k$-th plateau, and the integer interval $T_k:=[t_k,e_k]\cap\mathbb N_0$ the $k$-th transition. Their lengths are $N_k$ and $L_k$, respectively. The formulas defining $b_j$ agree at the common endpoints.

\begin{figure}[h]
    \centering 
    \includegraphics[width=0.7\textwidth]{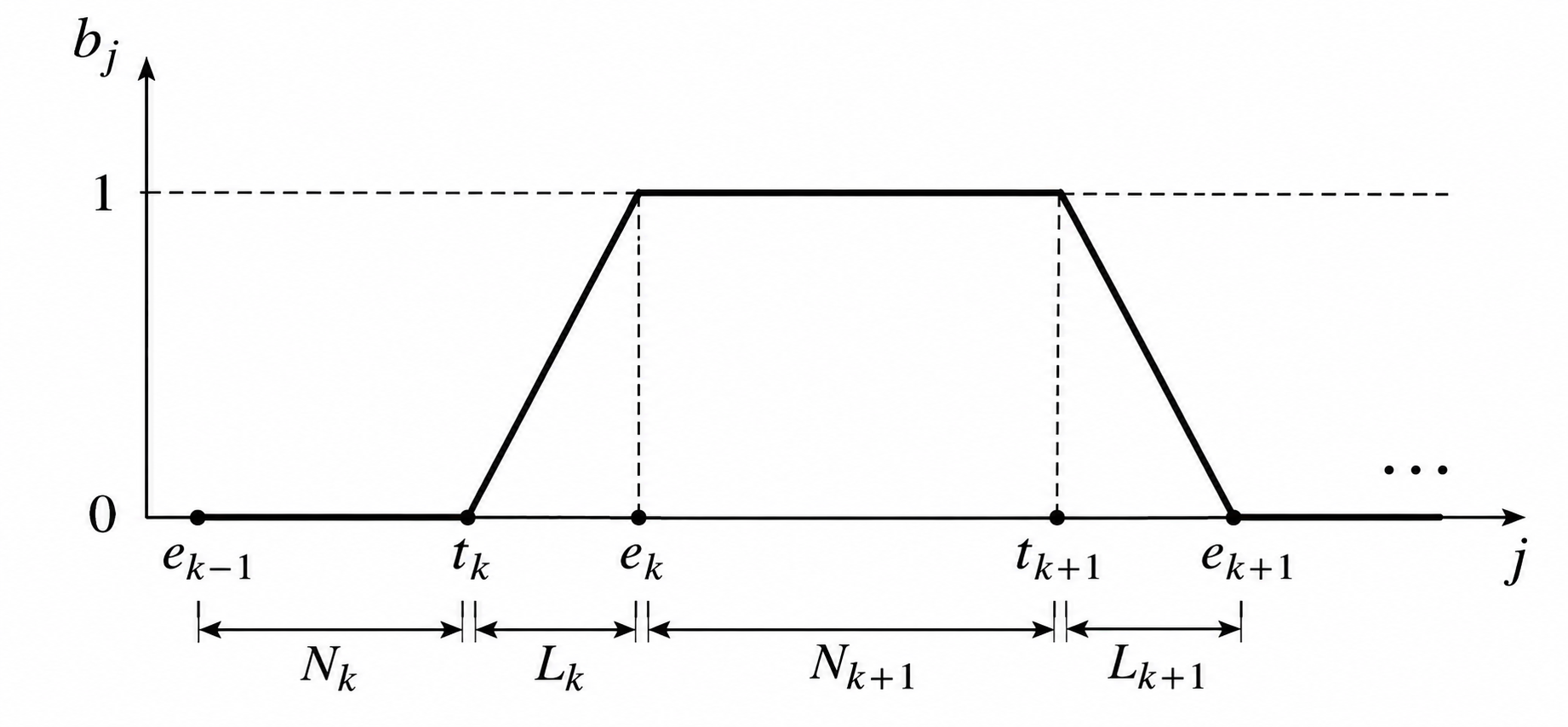}
\end{figure}
These formulas agree at boundaries. Moreover, for $e_{k-1}\leqslant j<e_k$, we have:
\[
|b_{j+1}-b_j|\leqslant\frac{1}{L_k}
\]
Since $L_k\to\infty$, it follows that $|b_{j+1}-b_j|\to 0$, $j\to\infty$.
Put $z^j=(b_j,b_{j+1},\ldots)$ and $Z=\overline{\{z^j:j\geqslant0\}}$.

To obtain the lower bound, we consider histories starting at different positions on the same plateau. We choose the initial positions at distance $L_k$
from one another. When one of these histories reaches the end of the transition, any history starting
further back is still on the plateau. Thus, the two histories differ by $1$ at some coordinate. In this way, we obtain approximately $N_k/L_k$ pairwise separated histories.

\begin{lemma}\label{lem:lower}
$\hp(\sigma|_{Z})\geqslant\alpha$.
\end{lemma}

\begin{proof}
For each $k\geqslant 1$, put
\[
m_k=\left\lfloor\frac{N_k}{L_k}\right\rfloor,
\qquad
n_k=N_k+L_k+1.
\]
We choose $m_k$ starting positions on the $k$th
plateau, spaced $L_k$ apart, by setting
\[
j_i=t_k-iL_k,\qquad 1\leqslant i\leqslant m_k.
\]
Since $m_kL_k\leqslant N_k$, we have
\[
e_{k-1}=t_k-N_k\leqslant j_i\leqslant t_k.
\]
Thus, all the chosen starting positions belong to
the $k$th plateau.

We claim that the set
\[
E_k=\{z^{j_i}:1\leqslant i\leqslant m_k\}
\]
is $(n_k,1/2)$-separated with respect to $\rho_{n_k}$.

Take any $1\leqslant i<\ell\leqslant m_k$. Since $j_i>j_\ell$,
the history starting at $j_i$ is closer to the end
of the plateau and therefore reaches the next
plateau earlier. We choose the observation time $r$ at which the
history starting at $j_i$ reaches the end of the
transition, that is,
\[
r=e_k-j_i
  =(t_k+L_k)-(t_k-iL_k)
  =(i+1)L_k.
\]
At this time, we have $j_i+r=e_k$ and hence $b_{j_i+r}=c_{k+1}.$

On the other hand, $j_\ell+r
=t_k-\ell L_k+(i+1)L_k
=t_k-(\ell-i-1)L_k.$ Since $\ell\geqslant i+1$, it follows that $e_{k-1}\leqslant j_\ell\leqslant j_\ell+r\leqslant t_k.$
Thus, at time $r$, the history starting at $j_\ell$
is still on the $k$th plateau, and consequently $b_{j_\ell+r}=c_k.$

As the values on two consecutive plateaus differ
by $1$, we obtain $|b_{j_i+r}-b_{j_\ell+r}|
=|c_{k+1}-c_k|=1.$ We only have to check if $r$ lies within the
observation window of length $n_k$. Since $i+1\leqslant\ell\leqslant m_k$, we have $0\leqslant r=(i+1)L_k
\leqslant m_kL_k\leqslant N_k<n_k.$
Therefore
\[
\rho_{n_k}(z^{j_i},z^{j_\ell})
\geqslant |b_{j_i+r}-b_{j_\ell+r}|
=1>\frac12.
\]
Since $i$ and $\ell$ were chosen arbitrarily,
$E_k$ is $(n_k,1/2)$-separated. Hence
\[
\operatorname{sep}
\left(Z,\rho_{n_k},\frac12\right)
\geqslant m_k.
\]
Finally, by the choice of $N_k$ and $L_k$, we have
\[
m_k\sim N_k^\alpha
\qquad\text{and}\qquad
n_k\sim N_k.
\]
Consequently,
\[
\begin{aligned}
h_{\mathrm{pol}}(\sigma|_Z)
&\geq
\limsup_{k\to\infty}
\frac{\log\operatorname{sep}
\left(Z,\rho_{n_k},\frac12\right)}
{\log n_k}\\
&\geqslant
\limsup_{k\to\infty}
\frac{\log m_k}{\log n_k}
=\alpha.
\end{aligned}
\]
\end{proof}

To obtain the upper bound, we first observe that,
for a fixed observation time $n$, all sufficiently
late histories are close to constant ones.
It is therefore enough to consider finitely many
blocks. Within each block, histories starting
sufficiently far from the end of the plateau
are constant in their first $n$ coordinates.
Thus, we only need to estimate the number of
representatives for the histories starting near
the end of the plateau or on the transition itself.
Our aim is to show that, for each fixed $\varepsilon>0$,
the number of spanning histories grows at most as
$n^\alpha$, up to a logarithmic correction. As a direct consequence of the following lemma, we obtain that $\hp(\sigma|_{Z})\leqslant\alpha$.

\begin{lemma}\label{lem:upper}
For every $\varepsilon>0$ there exists a constant
$C_{\alpha,\varepsilon}>0$ such that
\[
\operatorname{span}(Z,\rho_n,\varepsilon)
\leqslant C_{\alpha,\varepsilon}
n^\alpha(1+\log n).
\]
\end{lemma}

\begin{proof}
Fix $0<\varepsilon<1$. Recall that
\[
N_k=2^{2^k},
\qquad
L_k=\lceil N_k^{1-\alpha}\rceil.
\]
Since the sequence $(L_k)$ is nondecreasing and
$b$ is constant on each plateau and linear on each
transition, we have $|b_{i+1}-b_i|\leqslant \frac{1}{L_k}$, for every $i\geqslant e_{k-1}$. Consequently, for any $e_{k-1}\leqslant j<j'$ and
$r\geqslant 0$, we obtain
\[
|b_{j+r}-b_{j'+r}|
\leqslant
\sum_{i=j+r}^{j'+r-1}|b_{i+1}-b_i|
\leqslant \frac{j'-j}{L_k}.
\]
Taking the maximum over $0\leqslant r<n$, we get, for every $j,j'\geqslant e_{k-1}$:
\[
\rho_n(z^j,z^{j'})
\leqslant \frac{|j-j'|}{L_k}
\]
First, consider the blocks satisfying $L_k>\frac{4n}{\varepsilon}.$ For every starting index $j$ in such a block,
we have
\[
\begin{aligned}
\rho_n(z^j,(b_j,b_j,\ldots))
&=\max_{0\leq r<n}|b_{j+r}-b_j|\leqslant \frac{n-1}{L_k}
<\frac{\varepsilon}{4}.
\end{aligned}
\]
Thus, every such history is within
$\varepsilon/4$ of a constant sequence. Since $\rho_n((c,c,\ldots),(d,d,\ldots))=|c-d|,$ the family of all constant sequences with values
in $[0,1]$ can be covered by $O_\varepsilon(1)$
balls of radius $\varepsilon/4$.
By the triangle inequality, all histories starting
in blocks satisfying $L_k>4n/\varepsilon$ can
therefore be covered by $O_\varepsilon(1)$ balls
of radius $\varepsilon/2$. It remains to consider the blocks satisfying $L_k\leqslant \frac{4n}{\varepsilon}$.

Since $L_k\geqslant N_k^{1-\alpha}=2^{(1-\alpha)2^k},$ every such block satisfies
\[
2^{(1-\alpha)2^k}
\leq \frac{4n}{\varepsilon}.
\]
Hence
\[
k\leqslant
\log_2\left(
\frac{\log_2(4n/\varepsilon)}{1-\alpha}
\right).
\]
In particular, there are at most
$O_{\alpha,\varepsilon}(1+\log n)$ such blocks.

Consider one of these blocks. If the first $n$
coordinates of a history $z^j$ lie entirely
in the plateau, then $\rho_n(z^j,(c_k,c_k,\ldots))=0.$ Thus, all such histories can be covered by a
single ball centered at the constant sequence
$(c_k,c_k,\ldots)$. The remaining histories start either in the
plateau sufficiently close to its end to reach
the transition within the first $n$ coordinates,
or in the transition itself.

There are at most $\min\{N_k,n\}$ starting
indices of the first type and at most $L_k$
starting indices of the second type.
Moreover, these indices form a consecutive
interval containing at most
\[
\min\{N_k,n\}+L_k
\]
indices.

Set
\[
s_k=\max\left\{
1,\left\lfloor\frac{\varepsilon L_k}{4}\right\rfloor
\right\}.
\]
Choose representatives among these starting
indices at intervals of length $s_k$, beginning
with the first index of the interval. If $s_k=1$, every starting index is chosen.
Otherwise, for every starting index $j$ there
exists a chosen index $j'$ such that
\[
|j-j'|<s_k
\leqslant \frac{\varepsilon L_k}{4}.
\]
By the Lipschitz estimate established above,
\[
\rho_n(z^j,z^{j'})
\leq \frac{|j-j'|}{L_k}
<\frac{\varepsilon}{4}.
\]
Thus, the chosen histories form an
$\varepsilon/2$-spanning family for the
histories under consideration.

Since $s_k\geqslant \frac{\varepsilon L_k}{8},$
the number of chosen representatives is at most
\[
\begin{aligned}
1+\frac{\min\{N_k,n\}+L_k}{s_k}
&\leqslant
1+\frac{8}{\varepsilon}
\frac{\min\{N_k,n\}+L_k}{L_k}\\
&=O_\varepsilon\left(
1+\frac{\min\{N_k,n\}}{L_k}
\right).
\end{aligned}
\]
We now estimate the last quotient.
If $N_k\leqslant n$, then
\[
\frac{\min\{N_k,n\}}{L_k}
=\frac{N_k}{L_k}
\leqslant \frac{N_k}{N_k^{1-\alpha}}
=N_k^\alpha
\leqslant n^\alpha.
\]
On the other hand, if $N_k>n$, then
\[
\frac{\min\{N_k,n\}}{L_k}
=\frac{n}{L_k}
\leqslant \frac{n}{N_k^{1-\alpha}}
\leqslant \frac{n}{n^{1-\alpha}}
=n^\alpha.
\]
Consequently, in both cases,
\[
\frac{\min\{N_k,n\}}{L_k}\leq n^\alpha.
\]

Therefore, the histories starting in each of
the remaining blocks can be covered by
$O_\varepsilon(n^\alpha)$ balls of radius
$\varepsilon/2$. Combining this estimate with the bound on the
number of remaining blocks and the cover of
all histories starting in the later blocks,
we obtain
\[
\operatorname{span}
\left(\{z^j:j\geqslant 0\},\rho_n,\frac{\varepsilon}{2}\right)
\leqslant C_{\alpha,\varepsilon}n^\alpha(1+\log n).
\]

Finally, since $Z=\overline{\{z^j:j\geqslant 0\}},$ enlarging the radii of the covering balls from
$\varepsilon/2$ to $\varepsilon$ gives a cover
of $Z$ with the same number of balls. Hence
\[
\operatorname{span}(Z,\rho_n,\varepsilon)
\leqslant C_{\alpha,\varepsilon}n^\alpha(1+\log n),
\]
which completes the proof.
\end{proof}

A point on a hair need not stay on hairs forever. It can reach the fixed base interval after finitely many iterates. Hence the symbolic space must also contain histories of the form $(b_j,\ldots,b_s,b_s,b_s,\ldots)$. In principle, choosing a stopping time could contribute an extra factor of $n$ and increase polynomial entropy by one. 

We now show that adding stopped histories does not
increase the polynomial entropy. The main observation
is that stopping a history at different positions
on the same plateau gives the same history.
Moreover, if two stopping positions belong to
the same monotone transition, the distance between
the resulting histories is controlled by the
difference between their terminal values.
Thus, it is enough to estimate the number of
transitions encountered during the first $n$ steps.

For each $t\geqslant 0$, define the stopping operator
$S_t\colon\Sigma\to\Sigma$ by
\[
(S_t(z))_r=z_{\min\{r,t\}},\qquad r\geqslant 0.
\]
Thus, $S_t(z)$ is the history obtained by
stopping $z$ at time $t$. In other words, $$S_t(z)=(z_0,\dots,z_{t-1},z_t,z_t,\dots).$$
Notice that, for every $n\geqslant 1$ and $t\geqslant 0$,
the stopping operator $S_t$ is 1-Lipschitz
with respect to $\rho_n$. Indeed,
\[
\begin{aligned}
\rho_n(S_t(z),S_t(w))
&=\max_{0\leqslant r<n}
|z_{\min\{r,t\}}-w_{\min\{r,t\}}|\\
&\leqslant \rho_n(z,w).
\end{aligned}
\]
For $n\geqslant 1$ and $\varepsilon>0$, define
\[
Q(n,\varepsilon)
:=
\sup_{z\in Z}
\operatorname{span}
\bigl(
\{S_t(z):0\leqslant t<n\},
\rho_n,\varepsilon
\bigr).
\]

\begin{lemma}\label{lem:Q}
For every $\varepsilon>0$ there exists a constant
$C_\varepsilon>0$ such that
\[
Q(n,\varepsilon)
\leqslant C_\varepsilon(1+\log n).
\]
\end{lemma}

\begin{proof}
Fix $\varepsilon>0$ and consider first an ordinary
history $z=z^j=(b_{j+r})_{r\geqslant 0}$.

If two stopping positions $j+t$ and $j+t'$
belong to the same plateau, then
\[
S_t(z)=S_{t'}(z).
\]
If they belong to the same monotone transition,
then
\[
\rho_n(S_t(z),S_{t'}(z))
\leqslant |b_{j+t}-b_{j+t'}|.
\]
Since all values of $b$ belong to $[0,1]$,
the stopped histories corresponding to a single
transition can be covered by $O_\varepsilon(1)$
balls of radius $\varepsilon/2$.

It remains to estimate the number of plateaus and
transitions intersecting the interval $[j,j+n-1]$. Recall that the length of the $k$-th plateau is
$N_k=2^{2^k}$. Since $N_k\leq n$ implies
$k\leq \log_2\log_2 n$, the total number of
plateaus of length at most $n$ is
$O(1+\log n)$. On the other hand, the interval $[j,j+n-1]$
can intersect at most two plateaus of length
greater than $n$. Indeed, if it intersected
three such plateaus, it would have to contain
the middle one entirely, which is impossible
since its length is greater than $n$.

Thus, the number of plateaus intersecting
$[j,j+n-1]$ is $O(1+\log n)$, uniformly in $j$. Moreover, since any two consecutive transitions
are separated by a plateau, the number of
transitions intersecting $[j,j+n-1]$ is at most
the number of intersected plateaus plus one. Consequently, the total number of plateaus and
transitions intersecting $[j,j+n-1]$ is
$O(1+\log n)$, uniformly in $j$.

Combining these estimates, we obtain
\[
\operatorname{span}
\bigl(
\{S_t(z^j):0\leqslant t<n\},
\rho_n,\varepsilon/2
\bigr)
\leqslant C_\varepsilon(1+\log n).
\]

Now let $z\in Z$. By the definition of $Z$,
there exists an ordinary history $z^j$ such that
\[
\rho_n(z,z^j)<\frac{\varepsilon}{2}.
\]
Since
\[
\rho_n(S_t(z),S_t(z^j))
\leqslant \rho_n(z,z^j)
\]
for every $0\leqslant t<n$, the triangle inequality
gives
\[
\operatorname{span}
\bigl(
\{S_t(z):0\leqslant t<n\},
\rho_n,\varepsilon
\bigr)
\leqslant C_\varepsilon(1+\log n).
\]

Taking the supremum over $z\in Z$ completes
the proof.
\end{proof}

We now use the estimates obtained in the previous
two lemmas to bound the spanning numbers of $Y$.
The idea is to first approximate an ordinary history
and then approximate its stopped version.
Since the stopping operators are 1-Lipschitz
with respect to $\rho_n$, these two approximations
can be combined. This gives an upper bound in terms
of the spanning numbers of $Z$ and the quantity
$Q(n,\varepsilon)$.
Note that the following lemma holds for an arbitrary
sequence $b$, independently of its particular
structure.
\begin{lemma}\label{lem:Qcover}
For any sequence $b$ in $[0,1]$ and all
$n\geqslant 1$, $\varepsilon>0$, we have
\[
\operatorname{span}(Y,\rho_n,4\varepsilon)
\leqslant
\operatorname{span}(Z,\rho_n,\varepsilon)
Q(n,\varepsilon).
\]
\end{lemma}

\begin{proof}
Let $m=\operatorname{span}(Z,\rho_n,\varepsilon).$
Choose a covering of $Z$ by $m$ balls of radius
$\varepsilon$ with centers in $\Sigma$.
From each ball intersecting $Z$, choose a
representative $z_i\in Z$.

Thus, for every $z\in Z$, there exists $i$
such that
\[
\rho_n(z,z_i)<2\varepsilon.
\]
Since the stopping operators are 1-Lipschitz
with respect to $\rho_n$, we have
\[
\rho_n(S_t(z),S_t(z_i))
\leqslant \rho_n(z,z_i)<2\varepsilon.
\]

By the definition of $Q(n,\varepsilon)$,
for each $i$ the family
\[
\{S_t(z_i):0\leqslant t<n\}
\]
can be covered by at most $Q(n,\varepsilon)$
balls of radius $\varepsilon$.

By the triangle inequality, the corresponding
balls of radius $3\varepsilon$ cover all histories
$S_t(z)$ with $z\in Z$ and $0\leqslant t<n$.

Moreover, for $t\geqslant n-1$, the histories
$S_t(z)$ and $S_{n-1}(z)$ have identical first
$n$ coordinates. The same is true of $z$ and
$S_{n-1}(z)$.

Hence all generating histories of $Y$ can be
covered by at most $mQ(n,\varepsilon)$ balls
of radius $3\varepsilon$.

Enlarging the radius to $4\varepsilon$, we obtain
a cover of their closure $Y$. Therefore,
\[
\operatorname{span}(Y,\rho_n,4\varepsilon)
\leqslant
\operatorname{span}(Z,\rho_n,\varepsilon)
Q(n,\varepsilon).
\]
\end{proof}

Combining the estimates for ordinary and stopped
histories, we obtain the desired polynomial entropy
for both symbolic spaces.

\begin{theorem}
\label{thm:symbolic}
For every $\alpha\in(0,1)$, the sequence $b$
constructed above satisfies
\[
\hp(\sigma|_{Z})=\hp(\sigma|_{Y})=\alpha.
\]
\end{theorem}

\begin{proof}
By Lemmas \ref{lem:lower} and \ref{lem:upper},
we have
\[
\hp(\sigma|_{Z})=\alpha.
\]

On the other hand, by Lemmas \ref{lem:Q}
and \ref{lem:Qcover}, together with
Lemma \ref{lem:upper}, we obtain
\[
\begin{aligned}
&\operatorname{span}(Y,\rho_n,4\varepsilon)\\
&\quad\leqslant
\operatorname{span}(Z,\rho_n,\varepsilon)
Q(n,\varepsilon)\\
&\quad\leqslant C_{\alpha,\varepsilon}
n^\alpha(1+\log n)^2.
\end{aligned}
\]
Since
\[
\lim_{n\to\infty}
\frac{\log(1+\log n)}{\log n}=0,
\]
it follows that
\[
\hp(\sigma|_{Y})\leqslant\alpha.
\]

The reverse inequality follows from the inclusion
$Z\subseteq Y$, since
\[
\hp(\sigma|_{Y})
\geqslant\hp(\sigma|_{Z})=\alpha.
\]
\end{proof}

\subsection{Main application}

The sequence $b$ constructed in the previous
sections has the desired polynomial entropy,
but it cannot be used directly in our dendrite
construction. Indeed, its values are not pairwise
distinct, whereas each term of the sequence $a$
determines the attachment point of a different hair.
Moreover, we require all attachment points to
belong to the interior of the base interval.

We therefore seek a sequence $a$ of pairwise
distinct points in $(0,1)$ obtained by slightly
perturbing the sequence $b$. The main purpose
of this section is to show that such a perturbation
does not change the polynomial entropy, provided
that
\[
|a_j-b_j|\longrightarrow 0.
\]

To prove this, we first show that removing
finitely many initial terms of a sequence does
not affect the polynomial entropy of its ordinary
and stopped history spaces.

We will now use a more precise notation, i.e. $X_a$, $Y_a$ and $Z_a$, to distinguish between these spaces induced by different sequences.

\begin{lemma}\label{lem:tail}
For any sequence $a$ and $j\geqslant0$, the ordinary and stopped history spaces of $a$ have the same polynomial entropies as those of the tail $a^{(j)}=(a_j,a_{j+1},\ldots)$.
\end{lemma}
\begin{proof}
The $j$-fold shift of $Z_a$ equals $Z_{a^{(j)}}$. Moreover, $\sigma^j(Y_a)=Y_{a^{(j)}}\cup\{(a_s,a_s,\ldots):s<j\}$. The added set is finite and fixed under the shift. Apply Lemmas \ref{lem:image} and \ref{lem:union}.
\end{proof}

\begin{lemma}\label{lem:perturb}
If $a_j,b_j\in[0,1]$ and $|a_j-b_j|\to 0$, then
\[
\hp(\sigma|_{Z_a})=\hp(\sigma|_{Z_b}),
\qquad
\hp(\sigma|_{Y_a})=\hp(\sigma|_{Y_b}).
\]
\end{lemma}

\begin{proof}
We prove the assertion for stopped histories.
The proof for ordinary histories is analogous.

Fix $\varepsilon>0$. Since $|a_i-b_i|\to 0$,
there exists $j\geqslant 0$ such that
\[
\sup_{i\geqslant j}|a_i-b_i|<\varepsilon.
\]
Consequently, for every $j\leqslant i\leqslant s\leqslant\infty$
and every $n\geqslant 1$, the corresponding generating
histories satisfy
\[
\rho_n(u_a^{i,s},u_b^{i,s})<\varepsilon.
\]

Let
\[
m=\operatorname{span}
(Y_{b^{(j)}},\rho_n,\varepsilon)
\]
and let $w_1,\ldots,w_m\in\Sigma$ be the centers
of a corresponding family of spanning balls.

For every generating history $u_a^{i,s}$ of
$Y_{a^{(j)}}$, there exists $k$ such that
\[
\rho_n(u_b^{i,s},w_k)<\varepsilon.
\]
By the triangle inequality,
\[
\begin{aligned}
\rho_n(u_a^{i,s},w_k)
&\leqslant \rho_n(u_a^{i,s},u_b^{i,s})
+\rho_n(u_b^{i,s},w_k)\\
&<2\varepsilon.
\end{aligned}
\]

Thus, the balls with the same centers and radius
$2\varepsilon$ cover all generating histories
of $Y_{a^{(j)}}$. Enlarging the radius to
$3\varepsilon$, we obtain
\[
\operatorname{span}
(Y_{a^{(j)}},\rho_n,3\varepsilon)
\leqslant
\operatorname{span}
(Y_{b^{(j)}},\rho_n,\varepsilon).
\]

Recall that
\[
\sigma^j(Y_a)=Y_{a^{(j)}}\cup P_j,
\]
where
\[
P_j=\{a_s^\infty:0\leqslant s<j\}.
\]
Using the same covering argument as in the proof
of Lemma \ref{lem:image}, now with respect to
the pseudometrics $\rho_n$, we obtain
\[
\begin{aligned}
&\operatorname{span}(Y_a,\rho_{n+j},12\varepsilon)\\
&\quad\leqslant C_{j,\varepsilon}
\operatorname{span}(\sigma^j(Y_a),\rho_n,3\varepsilon)\\
&\quad\leqslant C_{j,\varepsilon}
\bigl(
\operatorname{span}(Y_{a^{(j)}},\rho_n,3\varepsilon)
+j
\bigr)\\
&\quad\leqslant C_{j,\varepsilon}
\bigl(
\operatorname{span}(Y_{b^{(j)}},\rho_n,\varepsilon)
+j
\bigr)\\
&\quad\leqslant C_{j,\varepsilon}
\bigl(
\operatorname{span}(Y_b,\rho_n,\varepsilon)
+j
\bigr).
\end{aligned}
\]
Here, the last inequality follows from the
inclusion $Y_{b^{(j)}}\subseteq Y_b$.

Since $j$ and $C_{j,\varepsilon}$ are independent
of $n$, taking logarithms, dividing by $\log n$,
passing to the limit superior as $n\to\infty$,
and then letting $\varepsilon\to 0$, we obtain
\[
\hp(\sigma|_{Y_a})
\leqslant
\hp(\sigma|_{Y_b}).
\]
Interchanging the roles of $a$ and $b$ gives
the reverse inequality.

\end{proof}

\begin{theorem}\label{thm:main}
For every $\alpha\in(0,1)$, there exist a dendrite
$D_\alpha$ and a continuous map
$F_\alpha\colon D_\alpha\to D_\alpha$ such that
\[
\hp(F_\alpha)=\alpha.
\]
\end{theorem}

\begin{proof}
Fix $\alpha\in(0,1)$ and consider the sequence $b$
constructed in the previous sections. Choose inductively pairwise distinct points
$a_j\in(0,1)$ such that
\[
|a_j-b_j|<2^{-j-2}.
\]
This is possible since every neighborhood of
$b_j\in[0,1]$ contains infinitely many points
of $(0,1)$.

Clearly, $|a_j-b_j|\to0$. Moreover, since
$|b_{j+1}-b_j|\to0$, the inequality
$$|a_{j+1}-a_j|
\leqslant |a_{j+1}-b_{j+1}|+|b_{j+1}-b_j|+|b_j-a_j|
$$
implies that $|a_{j+1}-a_j|\to0$.

Let $h_j=2^{-j-1}$ and construct the dendrite
$D_a$ and the continuous map $F_a$ as in
Theorem \ref{thm:realization}.

By Theorem \ref{thm:realization},
Lemma \ref{lem:perturb}, and
Theorem \ref{thm:symbolic}, we obtain
\[
\hp(F_a)
=\hp(\sigma|_{Y_a})
=\hp(\sigma|_{Y_b})
=\alpha.
\]

Taking $D_\alpha=D_a$ and $F_\alpha=F_a$
completes the proof.
\end{proof}
\vspace{1cm}

The question of realizing values of polynomial entropy in the interval $(1,2)$ for continuous dendrite maps remains open. In our construction, the choice of the lengths of the transitions relies on the assumption that $\alpha<1$. Indeed, for $\alpha>1$, the choice $L_k=\lceil N_k^{1-\alpha}\rceil)$ gives $L_k=1$, and hence the differences between consecutive terms of the sequence $b$ no longer tend to zero. Thus, the construction of the sequence $b$ would have to be modified. It would be interesting to see whether the same dendrite construction, with a different choice of the sequence $a$, could also yield polynomial entropy between one and two.

\vspace{1cm}

\textbf{Declaration on the use of artificial intelligence}\\
During the preparation of this paper, the authors used ChatGPT (OpenAI) as a research assistant to explore mathematical ideas, develop and refine arguments, and improve the presentation of the results. The research process was guided by the authors, who critically evaluated the AI-generated suggestions. All mathematical statements and proofs were carefully reviewed and verified by the authors, who take full responsibility for the content of the paper.

\end{document}